\documentclass[a4paper,12pt]{article}
\usepackage[english]{babel}
\usepackage[T2A]{fontenc}
\usepackage[cp1251]{inputenc}
\usepackage[english]{babel}
\usepackage{amsthm}
\usepackage[tbtags]{amsmath}
\usepackage{amsfonts,amssymb}
\begin{document}

\newtheorem{theorem}{Theorem}
\newtheorem{lemma}{Lemma}
\newtheorem{proposition}{Proposition}
\newtheorem{Cor}{Corollary}

\begin{center}
{\large\bf Centrally Essential Factor Rings and\\ 
Subdirect Indecomposability}
\end{center}
\begin{center}
Oleg Lyubimtsev\footnote{Nizhny Novgorod State University, Nizhny Novgorod, Russia; email: oleg\_lyubimcev@mail.ru .},
Askar Tuganbaev\footnote{National Research University MPEI, Moscow, Russia; Lomonosov Moscow State University, Moscow, Russia; tuganbaev@gmail.com .}
\end{center}

\textbf{Abstract.} Let $R$ be a ring and let $J(R)$, $C(R)$ be its Jacobson radical and center, correspondingly. If $R$ is a centrally essential ring and the factor ring $R/J(R)$ is commutative, then any minimal right ideal is contained in the center $C(R)$. A right Artinian (or right Noetherian subdirectly indecomposable) centrally essential ring is a right and left Artinian local ring. We describe centrally essential Noetherian subdirectly indecomposable rings and centrally essential rings with subdirectly indecomposable center. We give examples of non-commutative subdirectly indecomposable, centrally essential rings.


The work of Oleg Lyubimtsev is supported by Ministry of Education and Science of the Russian Federation, project FSWR-2023-0034. The study of Askar Tuganbaev is supported by grant of Russian Science Foundation (=RSF), project 22-11-00052, https://rscf.ru/en/project/22-11-00052.

\textbf{MSC2020 database 16D25, 16R99}

\section{Introduction}\label{sec1}

We consider only associative unital non-zero rings. A ring $R$ is said to be \textsf{centrally essential} if either $R$ is commutative or for every non-central element $a$ of $R$, there exist non-zero central elements $x,y$ with $ax = y$. It is clear that a ring $R$ with center $C(R)$ is centrally essential if and only if the module $R_{C(R)}$ is an essential extension of the module $C(R)_{C(R)}$. In addition, any commutative ring is centrally essential. 

\textbf{Remark 1.1.} There exists a non-commutative ring $R$ such that all factor rings of $R$ are centrally essential (see Example 3.3). 

\textbf{Remark 1.2.} In \cite[Proposition 2.8]{MT20b}, it is proved that a right uniserial, right Artinian, centrally essential ring is left Artinian.

\textbf{Remark 1.3.} There exist right and left subdirectly indecomposable rings $R$ which are not centrally essential. For example, let $V$ be a vector space over a field $\mathbb{F}$ of characteristic $0$ or $p\neq 2$ with basis $\{e_1, e_2\}$, and let $\Lambda(V)$ be the exterior algebra of the space $V$, i.e., $\Lambda(V)$ is an algebra with operation $\wedge$, generators $e_1, e_2$ and defining relations $e_i\wedge e_j + e_j \wedge e_i = 0 \,\,\mbox {for all}\,\, i, j = 1, 2$. Then $R = \Lambda(V)$ is a non-commutative $4$-dimensional $\mathbb{F}$-algebra with basis $\{1, e_1, e_2, e_1 \wedge e_2\}$ and the center $C(R)$ of $R$ is $\{q_0\cdot 1 + q_1\cdot e_1\wedge e_2 \, \mid \, q_i\in \mathbb{F}\}$.
It is easy to see that $R$ is a right and left subdirectly indecomposable ring with subdirectly indecomposable center. In addition, the ring $R$ is not centrally essential; see \cite[Theorem 1]{MT19}.

In connection to Remarks 1.1, 1.2 and 1.3, we prove Theorems 1.4, 1.5 and 1.6 which are the main results of this paper.

\textbf{Theorem 1.4.} Let $R$ be a centrally essential ring.

\textbf{1.} If the factor ring $R/J(R)$ is commutative, then any minimal right ideal $S$ of the ring $R$ is contained in the center $C(R)$. In particular, all minimal right ideals of the ring $R$ are ideals and its right socle is contained in the center of the ring $R$.

\textbf{2.} If $R$ is a right Artinian ring (or right Noetherian subdirectly indecomposable) ring, then $R$ is a right and left Artinian local ring.

\textbf{3.} If the center of the ring $R$ is subdirectly indecomposable, then the ring $R$ is right and left subdirectly indecomposable.

\textbf{Theorem 1.5.} If $R$ is a right Noetherian ring such that all factor rings of $R$ are centrally essential, then $R$ is a subdirect product of local, right and left Artinian rings.

\textbf{Theorem 1.6.} Let $R$ be a right Noetherian ring with center $C(R) = C$. We set $\overline{C} = J(R)\backslash J(C)$. 

\textbf{1.} If the ring $R$ is subdirectly indecomposable, then the ring $R$ is centrally essential if and only if 
$\text{Ann}_{\overline{C}}J(C) = \text{Ann}_{\overline{C}}J(R) = 0$.

\textbf{2.} If the center $C$ of the ring $R$ is subdirectly indecomposable with core $H$, then the ring $R$ is centrally essential if and only if $rC\cap H\neq 0$ for any non-zero element $r\in R$.

We give some necessary definitions and notation.
 
For a ring $R$, the Jacobson radical of $R$ is denoted by $J(R)$. A module is said to be \textsf{subdirectly indecomposable} if it has the least non-zero submodule. A ring $R$ is said to be \textsf{right subdirectly indecomposable} (or \textsf{$RSI$-ring} \cite{{Desh71}}) if the intersection $H_r$ of all its non-zero right ideals is non-zero. The \textsf{left core} $H_{\ell}$ of the ring $R$ is defined similarly. A ring $R$ is said to be \textsf{subdirectly indecomposable} with \textsf{core} $H$ if the intersection $H$ of all its non-zero ideals is non-zero. A ring $R$ is said to be \textsf{centrally subdirectly indecomposable} if the center of $R$ is subdirectly indecomposable.

For a module $M$, the socle \textsf{$Soc\, M$} is the sum of all simple submodules in $M$; if $M$ does not contain simple submodules, then $Soc\, M = 0$, by definition. A module $M$ is said to be \textsf{uniserial} if any two submodules of $M$ are comparable with respect to inclusion.

For a ring $R$, the \textsf{Lie ring} $G(R)$ of $R$ consists of elements of the ring $R$ with multiplication $[x, y] = xy - yx$ for any $x, y\in R$. An element $r\in R$ is said to be \textsf{right  regular} (resp., \textsf{left regular}) if $rx\neq 0$ (resp.,  $xr\neq 0$) for every non-zero element $x\in R$. Right and left regular elements are called  \textsf{regular}. A non-zero element $r\in R$ is \textsf{$C(R)$-torsion free} if $rc\neq 0$ for every non-zero element $c\in C(R)$. If $S$ and $T$ are two subsets of the ring $R$, then $\text{Ann}_ST = \{s\in S\,|\, sT = 0\}$. If $\mathbb{F}$ is a field and $V$ is a vector $\mathbb{F}$-space with basis $e_1,\ldots e_n$, then the \textsf{exterior algebra} $\Lambda(V)$ is an algebra with operation $\wedge$, generators $e_1,\ldots e_n$, defining relations
$e_i \wedge e_j + e_j \wedge e_i = 0 \,\, \mbox {for all} \,\, i, j = 1,\ldots n$, and the basis 
$$
\{1, e_{i_1}\wedge e_{i_2}\wedge\ldots\wedge e_{i_k}\,\,\mid k = 1, 2,\ldots n;\,\, 1\le i_1 < i_2 <\ldots i_k\le n \}.
$$ 

\section{Proof of Main Results}\label{sec2}

\textbf{Proposition 2.1.} For a centrally essential ring $R$, the following conditions are equivalent.

\textbf{1)} The ring $R$ is right subdirectly indecomposable.

\textbf{2)} The ring $R$ is left subdirectly indecomposable.

\textbf{3)} The ring $R$ is subdirectly indecomposable.

In addition, the cores $H_r$, $H_{\ell}$ and $H$ coincide in each of cases \textbf{1)}--\textbf{3)}.

\textbf{Proof.} Let the ring $R$ be right subdirectly indecomposable and let $\{J_t\}_{t\in T}$ be the set of all non-zero left ideals of the ring $R$. Since the ring $R$ is centrally essential, each of the left ideals $J_t$ contains a non-zero central element $c_t$. Then 
$0\neq \cap_{t\in T}c_tR\subseteq \cap_{t\in T}J_t$. Therefore,  \textbf{1)} imples \textbf{2)}.
 
It is clear that \textbf{2)} implies \textbf{3)}. Since each of the ideals $c_tR$ is two-sided, it follows from \textbf{3)} that \textbf{1)} is true.

It is clear that $H_r = hR$ for any $0\neq h\in H_r$. Since $R$ is a centrally essential ring, we have that $0\neq hc = d$ for some central elements $c$ and $d$. Then $H_r = dR$ is a two-sided ideal. Therefore, $H_r = H$. The relation $H_{\ell} = H$ can be proved similarly.~$\square$

We remark that the analogue of Proposition 2.1 is not true for rings which are not centrally essential. For example, the ring of all $n\times n$ ($n\ge 2$) matrices over a division ring is subdirectly indecomposable but it is not right or left subdirectly indecomposable. In addition, there exists a right subdirectly indecomposable ring which is not left subdirectly indecomposable; see \cite[Example 4.5]{Desh71}.

It follows from Proposition 2.1 that it is sufficient to restrict ourself to centrally essential, subdirectly indecomposable rings in the study of right (or left) subdirectly indecomposable, centrally essential rings.

\textbf{Proposition 2.2.} If $R$ is a centrally essential ring and the factor ring $R/J(R)$ is commutative, then any minimal right ideal $S$ of the ring $R$ is contained in $C(R)$. In particular, all minimal right ideals of the ring $R$ are ideals and $Soc\, R_R$ is contained in the center of the ring $R$.

\textbf{Proof.} The simple $R$-module $S_R$ is a cyclic module with non-zero generator $s\in S$. We set $K = S\cap C(R)$. By assumption, there exist non-zero central elements $x,y\in C(R)$ such that $0\ne sx = y\in S$. Then $S = yR$ and $y\in K\neq 0$. Since $J(R)$ is the intersection of the annihilators of all simple right $R$-modules, we have that $SJ(R) = 0$. It follows from commutativity of the ring $R/J(R)$ that $ab - ba\in J(R)$ for any elements $a,b\in R$. In addition, $k(ab - ba) = 0$ for any $k\in K$. Since $k\in C(R)$, we have that $(ka)b = kba = b(ka)$. Therefore, $ka\in C(R)$. Further, it follows from $ka\in S$ that $K$ is a right ideal of the ring $R$. Since $S$ is a minimal right ideal, we have $S = K\subseteq C(R)$.~$\square$

\textbf{Proposition 2.3.} Let $R$ be a centrally essential ring. 

\textbf{1.} An element $r\in R$ is regular if and only if $r$ is $C(R)$-torsion free. 

\textbf{2, \cite[Lemma 2.2 (2.2.1]{MT20b}.} In the ring $R$, every one-sided zero-divisor is a two-sided zero-divisor.

\textbf{Proof.} \textbf{1.} It is clear that if $r$ is regular, then it is $C(R)$-torsion free. Conversly, let $rs = 0$ for some non-zero element $s\in R$. Then $s\notin C(R)$ and there exist $c, d\in C(R)$ such that $0\neq sc = d$. Since $r$ is $C(R)$-torsion free, we have that $rd\neq 0$ and $rsc\neq 0$. This is a contradiction.

\textbf{2.} If $rs = 0$ for some $0\neq s\in R$, then \textbf{1} implies that $rd = dr = 0$ for some $0\neq d\in C(R)$.~$\square$

\textbf{Proposition 2.4.} Let $R$ be a right Noetherian, centrally essential, subdirectly indecomposable ring with core $H$. 

\textbf{1.} $R$ is a right Artinian local ring.

\textbf{2.} The factor ring $R/J(R)$ is commutative.

\textbf{3.} The core $H$ is contained in the center $C(R)$.

\textbf{Proof.} \textbf{1.} It follows from Proposition 2.3(2) that any right regular element is left regular. Then the assertion follows from \cite[Theorem 2.1]{Desh71}.

\textbf{2.} The assertion follows from \textbf{1} and \cite[Proposition 2.4]{MT19b}.

\textbf{3.} By \textbf{2}, the ring $R/J(R)$ is commutative. Now we use Proposition 2.2.~$\square$

\textbf{Proposition 2.5.} Let $R$ be a centrally essential ring. 

\textbf{1.} If the ring $R$ is right Artinian, then $R$ is left Artinian.

\textbf{2.} If $R$ is a right Noetherian, subdirectly indecomposable ring, then $R$ is a right and left Artinian ring.

\textbf{3.} If the ring $R$ is local and $rJ(C) = 0$ for $r\in R$, then $r\in J(C)$. 

\textbf{4.} If the center of the ring $R$ is subdirectly indecomposable, then $R$ is right and left subdirectly indecomposable.

\textbf{Proof.} \textbf{1.} It follows from \cite[Proposition 2.4]{MT19b} that the ring $R/J(R)$ is commutative. Then $[x, y]\in J(R)$ for any $x, y\in G(R)$. In addition, the radical $J(R)$ is nilpotent with some nilpotence index $n$. Then
$$
[x_1, y_1]\cdot\ldots\cdot [x_n, y_n] = 0
$$
for all elements $[x_i, y_i]\in G(R)$, $i = 1, \ldots n$. In this case, the Lie ring $G(R)$ is nilpotent. Then the result follows from \cite[Theorem 2.2]{Bjork71}.

\textbf{2.} It follows from Proposition 2.4(1) that the ring $R$ is right Artinian. By \textbf{1}, the ring $R$ is left Artinian.

\textbf{3.} We set $C = C(R)$. We assume that $rJ(C) = 0$. Since $R$ is a centrally essential ring, we have $0\neq rc = d\in J(C)$. 
It follows from the relation $J(C) = C\cap J(R)$ that the element $c$ is invertible. Then $r = c^{-1}d\in J(C)$.

\textbf{4.} Let $M = \{I_k\}_{k\in K}$ be the set of all non-zero right ideals of the ring $R$. Since $R$ is a centrally essential ring, each of the ideals $I_k$ contains a non-zero central element $c_k$. We consider ideals $c_kC(R)$ of the center $C(R)$. Since the center is subdirectly indecomposable, we have
$0\neq \cap_{k\in K}c_kC(R)\subseteq \cap_{k\in K}I_k$.~$\square$

\textbf{2.6. The proof of Theorem 1.4.} 

\textbf{1.} The assertion follows from Proposition 2.2.

\textbf{2.} The assertion follows from Proposition 2.4(1) and Proposition 2.5(1,2).

\textbf{3.} The assertion follows from Proposition 2.5(4).~$\square$

\textbf{2.7. The proof of Theorem 1.5.} 

Any ring is a subdirect product of subdirectly indecomposable rings. Therefore, the right Noetherian ring $R$ is a subdirect product of right Noetherian, subdirectly indecomposable, centrally essential rings $R_i$, $i\in I$. It follows from Propositions 2.4(1) and 2.5(2) that all these rings are right and left Artinian local rings.~$\square$

\label{!!!}

\textbf{2.8. The proof of Theorem 1.6.} 
We recall that $R$ is a right Noetherian ring with center $C$ and $\overline{C} = J(R)\backslash J(C)$. 

\textbf{1.} If $R$ is a centrally essential ring, then it follows from Propositions 2.4(1) and 2.5(3) that $\text{Ann}_{\overline{C}}J(C) = 0$. Then $\text{Ann}_{\overline{C}}J(R) = 0$.

Conversely, let $0\neq r\notin C$ and $rJ(C)\neq 0$. By Proposition 2.4(1), the ring $R$ is local; therefore, $r\in\overline{C}$ or $r$ is invertible. Then $R$ is right Artinian and we can choose central elements $c_1,\ldots, c_{k-1}$ such that $rc_1\ldots c_{k-1}\neq 0$, but $rc_1\ldots c_{k-1}c = 0$ for all $c\in J(C)$. We set $c' = c_1\ldots c_{k-1}$. If $0\neq rc'\notin J(C)$, then $rc'\in \text{Ann}_{\overline{C}}J(C) = 0$. This is a contradiction. Consequently, $rc'$ is a non-zero central element and $R$ is a centrally essential ring.

\textbf{2.} Let $R$ be a centrally essential ring and $r\notin H$. There exist central elements $c$ and $d$ such that $0\neq rc = d$. If $d\in H$, then it is nothing to prove. Let $d\notin H$. Since $R$ is a centrally subdirectly indecomposable ring, we have 
$dC\cap H\neq 0$, i.e., $0\neq dc' = h$ for some $c'\in C$ and $h\in H$. Then
$$
0\neq dc' = (rc)c' = r(cc') = rc''\in rC,
$$
whence we obtain that $0\neq rc''\in rC\cap H$.

Conversely, for any $0\neq r\in R$, there exists an element $c\in C$ such that $0\neq rc = h\in H\subseteq C$, i.e., $R$ is a centrally essential ring.~$\square$

\section{Examples and Remarks}\label{sec3}

A ring is said to be \textbf{invariant} if all its one-sided ideals are two-sided ideals.

\textbf{Example 3.1.} We give an example of an Artinian, centrally essential, right and left subdirectly indecomposable, invariant ring $R$ such that $R$ is not a centrally subdirectly indecomposable ring and all proper centrally essential factor rings of $R$ are commutative. However,, there exists a proper factor ring of the ring $R$ which is not subdirectly indecomposable or centrally essential. A ring $R$ is a non-commutative $8$-dimensional algebra over an arbitrary field $\mathbb{F}$ of characteristic $0$ or $p\neq 2$. If we take $\mathbb{F} = \mathbb{Z}/p\mathbb{Z}$ for prime integer $p > 2$, then the ring $R$ is finite.

\textbf{Proof.} Let $R = \Lambda(V)$ be the the exterior algebra with generators $e_1, e_2, e_3$. Then $R$ is an $8$-dimensional $\mathbb{F}$-algebra with basis 
$$
\{1, e_1, e_2, e_3, e_1 \wedge e_2, e_2 \wedge e_3, e_1 \wedge e_3, e_1 \wedge e_2 \wedge e_3\}
$$ 
and $R$ is a non-commutative centrally essential ring (see details in \cite[Example 1]{MT19}) with center
$$
C(R) = \{q_0\cdot 1 + q_1\cdot e_1\wedge e_2 + q_2\cdot e_2\wedge e_3 + q_3\cdot e_1\wedge e_3 + q_4\cdot e_1 \wedge e_2 \wedge e_3, \, q_i\in \mathbb{F}\}.
$$
However, $I_1\cap I_2 = \{0\}$ for central ideals $I_1 = \langle e_1\wedge e_2\rangle_\mathbb{F}$ and 
$I_2 = \langle e_1\wedge e_2\wedge e_3\rangle_\mathbb{F}$, i.e., the center of the ring $R$ is not a subdirectly indecomposable ring. We remark that $R$ is an invariant ring and the ideal $H = I_2$ is contained in the intersection of any ideals of $R$. Consequently, $R$ is an Artinian, centrally essential, right and left subdirectly indecomposable ring which is not centrally subdirectly indecomposable.

We prove that all proper centrally essential factor rings of $R$ are commutative. We consider the factor ring $\overline{R} = R/I$, where 
$I$ is a non-zero ideal in $R$. If $e_i\notin I$ and $\overline{e_i}\notin C(\overline{R})$ for some $i$, then 
$\overline{e_i} J(C(\overline{R})) = \overline{0}$, since $I$ contains the core $H$ and $e_i\wedge (e_i\wedge e_j) = 0$. Then the ring $\overline{R}$ is not centrally essential, by Proposition 2.5(3). In addition, if, for example, 
$I = \langle \overline{e_1}\wedge\overline{e_3}, \overline{e_1}\wedge\overline{e_2}\wedge\overline{e_3}\rangle_\mathbb{F}$, then the ring $\overline{R}$ is not subdirectly indecomposable. Indeed, principal ideals generated by elements $\overline{e_1}\wedge\overline{e_2}$ and $\overline{e_2}\wedge\overline{e_3}$ have the zero intersection.~$\square$

\textbf{Example 3.2.} Let $Q_8$ be the quaternion group, i.e., the group with two generators $a$, $b$ and defining relations 
$a^4 = 1$, $a^2 = b^2$ and $aba^{-1} = b^{-1}$; e.g., see \cite[Section 4.4]{H59}. We have 
$$
Q_8 = \{e, a, a^2, b, ab, a^3, a^2b, a^3b\},
$$
with center $C(Q_8) = \{e, a^2\}$. We consider the group algebra $R = \mathbb{Z}_2 Q_8$ with center
$$
C(R) = \{\lambda_0 + \lambda_1a^2 + \lambda_2(a + a^3) + \lambda_3(b + a^3b) + \lambda_4(ab + a^3b),\quad \lambda_i\in \mathbb{Z}_2\};
$$
e.g., see \cite[Part 2]{P77}. Then $R$ is a ring with subdirectly indecomposable center and core $H = \{0, \widehat Q_8\}$, where
$\widehat Q_8 = \sum_{g\in Q_8}g$. In addition, $R$ is a centrally essential ring; see, \cite{MT18}. It follows from Theorem 1.6(2) that for every $0\neq x\in R$, there exists an element $c\in C(R)$ such that $xc = \widehat Q_8$. Then the factor ring $\overline{R} = R/H$ with basis 
$\{\overline{e}, \overline{a}, \overline{a}^2, \overline{a}^3, \overline{b}, \overline{ab}, \overline{a^2b}\}$ and with center 
$$
C(\overline{R}) = \{\alpha_0\overline{e} + \alpha_1\overline{a}^2 + \alpha_2(\overline{a} + \overline{a}^3) + 
\alpha_3(\overline{b} + \overline{a^2b}) + \alpha_4(\overline{e} + \overline{a} + \overline{a}^2 + \overline{a}^3 + \overline{b} + \overline{a^2b}\}
$$
is not a centrally essential ring. This property can be directly verified. For example, let $x\in R$ and $xa\neq ax$ but 
$\overline{x}\, \overline{a} = \overline{a}\, \overline{x}$ in $\overline{R}$. Then $xa = ax + \widehat Q_8$, whence $\sum x_ia = \sum ax_i + \widehat Q_8$. Since only basis elements $b$, $ab$, $a^2b$ and $a^3b$ do not commute with the element $a$, we obtain the contradictory relation
$$
\alpha_1b + \alpha_2ab + \alpha_3a^2b + \alpha_4a^3b =\widehat Q_8,\quad \alpha_i\in \mathbb{Z}_2.
$$
The case $xb\neq bx$ is considered similarly, but $\overline{x}\, \overline{b} = \overline{b}\, \overline{x}$. When we multiply the element 
$\overline{r} = (\overline{e} + \overline{a} + \overline{b} + \overline{ab})\notin C(\overline{R})$ by central basis elements, we obtain either $\overline{r}$ (under multiplication by $\overline{a}^2$) or $\overline{0}$. Consequently, $\overline{R}$ is not a centrally essential ring. Similarly, we can verify that all centrally essential proper factor rings of the ring $R$ are commutative.

\textbf{Example 3.3.} We give an example of a non-commutative ring $R$ with subdirectly indecomposable center such that all factor rings of $R$ are centrally essential.
Let $F = \mathbb{Q}(x,y)$ be the field of rational functions. We consider two partisl derivations 
$d_1 = \dfrac{\partial}{\partial x}$ and $d_2=\dfrac{\partial}{\partial x}$. Then the ring $R = T(F,F)$ of matrices
$$
\left\lbrace\left.\begin{pmatrix}
f&d_1(f)&g\\
0&f&d_2(f)\\
0&0&f
\end{pmatrix}\;\;\right|\;\; f,g\in F \right\rbrace
$$
is a centrally essential ring; see \cite{MT19b}. In addition, $R$ is a ring with subdirectly indecomposable center and with core
$$
H = \left\lbrace\left.\begin{pmatrix}
0&0&g\\
0&0&0\\
0&0&0
\end{pmatrix}\;\;\right|\;\; g\in F \right\rbrace.
$$
Since $R/H \cong \mathbb{Q}(x,y)$ and $H$ is the the least ideal in $R$, the ring $R$ does not have non-zero proper ideals not equal to $H$. 

\textbf{Remark 3.4.} In \cite[Proposition 3.3]{MT19}, it is proved that if $R$ is a centrally essential local ring, then the ring $R/J(R)$ is commutative (in particular, $R/J(R)$ is centrally essential). This is not true in general case; see \cite[Theorem 1.5]{MT20c}. 

\textbf{Remark 3.5.} If $R$ is not centrally essential, then the asserton of Proposition 2.4(1) is not true: there exist right Noetherian, right subdirectly indecomposable rings which are not right Artinian; see \cite[Example 2.2]{Desh71}. 

\textbf{Remark 3.6.} If $R$ is a centrally essential, right and left subdirectly indecomposable ring with core $H$, then either $H^2 = 0$ or $R$ is a field. Indeed, if $H^2 = H$, then $R$ is a division ring; see \cite[Proposition 1.3]{Desh71}. Since $R$ is a centrally essential ring, it follows from \cite[Proposition 3.3]{MT18} that $R$ is a field.

\textbf{Open question 3.7.} Is it true that there exists a centrally essential ring which has a proper non-commutative centrally essential factor ring?

\textbf{Open question 3.8.} Is it true that a centrally essential subdirectly indecomposable (or centrally subdirectly indecomposable) ring is invariant?

\label{!!!}

\end{document}